\documentclass[11pt]{article}
\usepackage[english]{babel}
\usepackage[T1]{fontenc}
\usepackage{verbatim}
\usepackage{amsmath}
\usepackage{amssymb}
\usepackage{graphicx}

\newcommand{\R}{\mathbb R}

\newcommand{\C}{   {\ifmmode{{\mathbb C}}\else{$\mathbb C$}\fi}}

\newcommand{\norm}[1]{\left\Vert #1\right\Vert}

\newcommand{\Z}{\mathcal{Z}}

\newcommand{\n}{\mathfrak{N}}

\newcommand{\X}{\mathcal{X}}
\newcommand{\op}{\overrightarrow}

\newcommand{\g}{\mathfrak{g}}
\newcommand{\z}{\mathfrak{z}}

\newcommand{\dt}{\frac{d}{dt}}
\newcommand{\dto}{\frac{d}{dt}_{|t=0}}
\newcommand{\dso}{\frac{d}{ds}_{|s=0}}

\newcommand{\Span}{\mbox{Span}}
\newcommand{\ad}{\mbox{ad}}

\newcommand{\ds}{\displaystyle}
\newcommand{\demo}{\noindent\textit{Proof. }}

\newcommand{\sig}{$(\Sigma)$ }
\newcommand{\sigp}{$(\Sigma')$ }

\begin{document}

\title{\bf Isometries of almost-Riemannian structures on Lie groups}

\author{Philippe JOUAN\footnote{Lab. R.~Salem, CNRS UMR 6085, Universit\'e
    de Rouen, avenue de l'universit\'e BP 12, 76801
    Saint-\'Etienne-du-Rouvray France. E-mail: Philippe.Jouan@univ-rouen.fr}, Guilherme ZSIGMOND\footnote{University of Rouen and Universidad Cat\'olica del Norte, Antofagasta Chile}, Victor AYALA\footnote{Instituto de Alta Investigaci\'on, Universidad de Tarapac\'a, Arica, Chile. E-mail: vayala@ucn.cl}}

\date{\today}

\maketitle
\begin{abstract}
A simple Almost-Riemannian Structure on a Lie group $G$ is defined by a linear vector field (that is an infinitesimal automorphism) and $\dim(G)-1$ left-invariant ones. 

It is first proven that two different ARSs are isometric if and only if there exists an isometry between them that fixes the identity. Such an isometry preserves the left-invariant distribution and the linear field. If the Lie group is nilpotent it is an automorphism.

These results are used to state a complete classification of the ARSs on the 2D affine and the Heisenberg groups.

\vskip 0.2cm


\vskip 0.2cm
\end{abstract}


\section{Introduction}

This paper is devoted to isometries of Almost-Riemannian structures on Lie groups. The purpose is to classify these structures, to find geometric invariants, and to determine their groups of isometries.

An almost-Riemannian structure (ARS in short) on an $n$-dimensional differential manifold can be defined, at least locally, by a set of $n$ vector fields, considered as an orthonormal frame, that degenerates on some singular set. This geometry goes back to \cite{Grushin} and \cite{Takasu}. It appears as a part of sub-Riemannian geometry, and has aroused some interest, as shown by the recent papers \cite{ABS08}, \cite{ABCGS10}, \cite{BCST09}, \cite{BCGJ11}, \cite{BCGM14}, \cite{BCGS}, \cite{BCG13}.

On an $n$-dimensional connected Lie group the simplest ARSs are defined by a set of $n-1$ left-invariant vector fields and one linear vector field, the rank of which is equal to $n$ on a proper open and dense subset and that satisfy the rank condition (a vector field on a Lie group is linear if its flow is a one parameter group of automorphisms, see Section \ref{BD}, and \cite{AS01}, \cite{AT99},  \cite{DJ14}, \cite{Jou11}, \cite{Jou12} about linear systems on Lie groups).

These ARSs, among which we find the famous Grushin plane  on the Abelian Lie group $\R^2$, has been studied in \cite{AJ16}. Among the results of this paper there is a study of the singular locus, that is the set of points where the vector fields fail to be independent. It is an analytic set, but not a subgroup, not even a submanifold, in general. However sufficient conditions for the singular locus to be a submanifold or a subgroup were exhibited and are recalled in Section \ref{BD}. This locus is very important in what concern the structure of ARSs, in particular in view of a classification. Another important geometric locus is the set of singularities of the linear field. It is always a subgroup.

In this paper we deal with smooth isometries, i.e. diffeomorphisms that respect the Euclidean metric of the tangent space at each point. First of all we show that such an isometry should preserve the singular locus and the group of singularities of the linear field. The main consequence is that the group of isometries of an ARS does not act transitively on $G$. Another consequence is that a left translation $L_g$ is an isometry if and only if $g$ belongs to the set of singularities of the linear field.

Then we prove that the isometries preserve the left-invariant distribution generated by the $n-1$ left-invariant vector fields, and also the linear field (up to the sign), see Theorem \ref{preservation}.

These constraints are rather strong. Consider for instance ARSs on the Heisenberg group, in the case where the distribution generated by the left-invariant vector fields is a subalgebra. Then the group of isometries is  generically reduced to the identity.

Kivioja and Le Donne proved in \cite{KL16} that the isometries of left-invariant metrics on nilpotent groups are affine, that is composed of an automorphism and a left translation. Their result cannot be directly applied here since the metric is not left-invariant. It can however be adapted, and it is shown in Theorem \ref{Nilpotent} that the isometries of ARSs on nilpotent Lie groups are affine.

The paper is organized as follows. In Section \ref{BD} the basic definitions and notations are stated, together with some useful results.

Section \ref{General} is devoted to the general theorems quoted above.

In Sections \ref{Aff} and \ref{Heisenberg} the ARSs of the 2D affine group and the 3D Heisenberg group are completely classified. We begin by a classification by isometries, up to a global rescaling of the metric that does not modify the geometry.
In order to reduce the number of parameters we then identify ARSs that differ only by the left-invariant metric on the left-invariant distribution. The typical ARSs that we obtain are characterized by the following geometric invariants: the singular locus, the set of singularities of the linear field, the tangency points and the eigenvalues of the derivation associated to the linear field.



\section{Basic definitions}\label{BD}
\subsection{Linear vector fields}

In this section the definition of linear vector fields and some of their properties are recalled. More details can found in \cite{Jou09}.

Let $G$ be 
a connected Lie group and $\g$ its Lie algebra (the set of left-invariant vector fields, identified with the tangent space at the identity). A vector field on $G$ is said to be {\it linear if its flow is a one-parameter group of automorphisms}. Notice that a linear vector field is consequently analytic and complete.
The flow of a linear vector field $\X$ will be denoted by $(\varphi_t)_{t\in\R}$.

The following characterization will be useful in the sequel.

\vskip 0.2cm
\noindent \textit{\textbf{Characterization of linear vector fields}}

{\it A vector field $\X$ on a connected Lie group $G$ is linear if and only if $\X$ belongs to the normalizer of $\g$ in the algebra $V^{\omega}(G)$ of analytic vector fields of $G$ and verifies $\X(e)=0$,
that is
\begin{equation} \label{deriv}
\forall Y\in \g \qquad [\X,Y]\in\g \quad\mbox{ and }\quad \X(e)=0.
\end{equation}
}

According to (\ref{deriv}) one can associate to a given linear vector field $\X$ the derivation $D$ of $\g$ defined by:
$$
\forall Y\in \g \qquad DY=-[\X,Y],
$$
that is $D=-\ad(\X)$.
The minus sign in this definition comes from the formula
$[Ax,b]=-Ab$ in $\R^n$. It also enables to avoid a minus sign in the useful formula:
\begin{equation}\label{equationfonda}
\forall Y\in\g, \quad \forall t\in \R \qquad \varphi_t(\exp
Y)=\exp(e^{tD}Y).
\end{equation}

An \textit{affine vector field} is an element of the normalizer $\n$ of $\g$ in $V^{\omega}(G)$, that is
$$
\n=\mbox{norm}_{V^{\omega}(G)}\g=\{F\in V^{\omega}(G);\ \forall Y\in \g, \quad
[F,Y]\in\g\},
$$
so that an affine vector field is linear if and only if it vanishes at the identity.

It can be shown (see \cite{AT99} or \cite{Jou09}) that an affine vector field can be uniquely decomposed
into a sum $F=\X+Z$
where $\X$ is linear and $Z$ is right-invariant.


\subsection{Almost-Riemannian structures}
For all that concern general sub-Riemannian geometry, including almost-Riemannian one, the reader is referred to \cite{ABB14}. About Almost-Riemannian structures on Lie groups more details can be found in \cite{AJ16}.

\newtheorem{DefARS}{Definition}
\begin{DefARS} \label{DefARS}
An almost-Riemannian structure on an $n$-dimensional Lie group $G$ is defined by a set of $n$ vector fields $\{\X,Y_1,\dots,Y_{n-1}\}$ where
\begin{enumerate}
	\item[(i)] $\X$ is linear;
	\item[(ii)] $Y_1,\dots,Y_{n-1}$ are left-invariant;
	\item[(iii)] $n=\dim G$ and the rank of $\X,Y_1,\dots,Y_{n-1}$ is full on a nonempty subset of $G$;
	\item[(iv)] the set $\{\X,Y_1,\dots,Y_{n-1}\}$ satisfies the rank condition.
\end{enumerate}

The metric is defined by declaring the frame $\{\X,Y_1,\dots,Y_{n-1}\}$ to be orthonormal.
\end{DefARS}

\vskip 0.2cm

\noindent \textbf{Equivalent definition}. \textit{An ARS on a Lie group can as well be defined by an $(n-1)$-dimensional left-invariant distribution $\Delta$ (that is $\Delta=\Span\{Y_1,\dots,Y_{n-1}\}$), a left-invariant Euclidean metric on $\Delta$ and a linear vector field $\X$ assumed to satisfy the conditions $(iii)$ and $(iv)$ of Definition \ref{DefARS}. The metric of the ARS is then defined by declaring $\X$ unitary and orthogonal to $\Delta$.}

\vskip 0.2cm

\noindent \textbf{Necessary conditions for the rank condition}

First notice that if $[\Delta,\Delta]\subseteq\Delta$ and $D(\Delta)\subseteq\Delta$, then the Lie algebra generated by $\X,Y_1,\dots,Y_{n-1}$ is  equal to $\R\X\oplus \Delta$. But the rank of that Lie algebra is not full at the identity $e$.
Consequently the rank condition implies that at least one of the following conditions hold:
\begin{enumerate}
	\item[(i)] $[\Delta,\Delta]\nsubseteq\Delta$ 
	\item[(ii)] $D(\Delta)\nsubseteq\Delta$
\end{enumerate}
In both cases, the full rank is obtained after one step.

\vskip 0.2cm

\noindent\textbf{Singular locus}

The set where the rank of $\X,Y_1,\dots,Y_{n-1}$ is not full will be referred to as the \textit{singular locus} and denoted by $\Z$. It is an analytic subset of $G$. By Assumption (iii) it is not equal to $G$, and by analycity its interior is empty. On the other hand $\X(e)=0$ and it cannot be empty. Finally $G\setminus \Z$ is an open, dense and proper subset of $G$.

The \textit{set of singularities of the linear field} $\X$ will be denoted by $\Z_\X=\{g\in G;\ \ \X(g)=0\}$. It is a subgroup of $G$, the Lie algebra of which is $\ker(D)$, where $D$ is the derivation associated to $\X$. It is included in $\Z$ but different in general.

The points of $G\setminus \Z$ will be called \textit{the Riemannian points} and the ones of $\Z$ \textit{the singular points}.

\vskip 0.2cm

\noindent\textbf{Norms and isometries}

The almost-Riemannian norm on $T_gG$ is defined by:
$$
\ds \mbox{For }\ X\in T_gG, \qquad\norm{X}=\min\left\{\sqrt{v^2+\sum_1^n u_i^2};\ v\X_g+u_1Y_1(g)+\dots+u_{n-1}Y_{n-1}(g)=X\right\}.
$$
It is infinite if the point $g$ belongs to the singular locus and $X$ does not belong to $\Delta_g$.

\newtheorem{Isometry}[DefARS]{Definition}
\begin{Isometry}\label{Isometry}
Let $(\Sigma)$ and $(\Sigma')$ be two ARSs on the Lie group $G$. An isometry $\Phi$ from $(\Sigma)$ onto $(\Sigma')$ is a diffeomorphism of $G$ that respects the norms, that is:
$$
\forall g\in G,\ \forall X\in T_gG \qquad \norm{T_g\Phi.X}_{\Sigma'}=\norm{X}_\Sigma.
$$
where $\norm{.}_\Sigma$ (resp. $\norm{.}_{\Sigma'}$) stands for the norm associated to $(\Sigma)$ in $T_gG$ (resp. to $(\Sigma')$ in $T_{\Phi(g)}G$).
\end{Isometry}


\subsection{Summary of some results}
The following results are proved in \cite{AJ16}.

\vskip 0.2cm

\noindent\textbf{Theorem 1 of \cite{AJ16}}
\textit{If $\Delta$ is a subalgebra of $\g$ then the singular locus $\Z$ is an analytic, embedded, codimension one submanifold of $G$.
Its tangent space at the identity is $D^{-1}\Delta$.
}

\vskip 0.2cm

\noindent\textbf{Theorem 2 of \cite{AJ16}}
\textit{If the Lie algebra $\g$ is solvable and $\Delta$ is a subalgebra of $\g$, then the singular locus $\Z$ is a codimension one subgroup of $G$ whose Lie algebra is $\z=D^{-1}\Delta$.}


\subsection{Notations}
In the sequel the following notations will be used:
\begin{enumerate}
\item Let $Y\in\g$. Then $Y_g$ stands for $TL_g.Y$, where $L_g$ is the left translation by $g$, and $TL_g$ its differential. The Euclidean norm in $T_gG$ is denoted by $\norm{.}_g$, or simply by $\norm{.}$ if no confusion is possible. Notice that if $g\in\Z$ then $\norm{.}_g$ is only defined on $\Delta_g$.
\item To the linear vector field $\X$ we associate $F(g)=TL_{g^{-1}}.\X_g\in\g$. In some expressions we also write $F_g$ for $F(g)$ when it is lighter, for instance in $\ad(F_g)$.
\end{enumerate}

Finally the following formula, proven in \cite{AJ16}, will be used in some proofs.
\begin{equation}\label{TgF}
T_gF=(D+\ad(F_g))\circ TL_{g^{-1}}.
\end{equation}
Notice that here $F_g$ is an element of $\g$.


\section{General Theorems about Isometries of ARSs}\label{General}

In this section \sig and \sigp stand for two simple ARSs on the same Lie group $G$. An object related to \sigp will be denoted by the same symbol as the analogous object related to \sig but with a prime, that is, $D'$, $(\varphi'_t)_{t\in\R}$, and so on.

\subsection{Translations}
The purpose of this subsection is to characterize the left translations that are isometries and to show that there exists an isometry from \sig onto \sigp  if and only if there exists one such isometry that preserves the identity.

\newtheorem{Prop1}{Proposition}
\begin{Prop1} \label{Prop1}
Let $\Phi$ be an isometry from \sig onto \sigp. Then $\Phi$ sends the singular locus of \sig onto the one of \sigp and the set of fixed points of $\X$ onto the one of ${\X'}$, that is:
$$
\Phi(\Z)=\Z' \qquad\mbox{and}\qquad \Phi(\Z_\X)=\Z'_{\X'}.
$$
In particular $\X'_{\Phi(e)}=0$.
\end{Prop1}

\demo
Since the singular locus is the set of points of $G$ where the rank of $\Span\{\X,Y_1,\dots, Y_{n-1}\}$ is not full, it is clear that $\Phi(\Z)=\Z'$. We are left to prove that $\Phi(\Z_\X)=\Z'_{\X'}$.

Let $V\in\Delta_e$, $V\neq 0$. It can be written in an unique way as $\ds V=\sum_{i=1}^{n-1}a_iY_i$ and we may assume without lost of generality that $\ds \norm{V}_e=\sqrt{\sum_{i=1}^{n-1}a_i^2}=1$.

Let $g\in \Z_\X$ and $V_g=TL_g.V$. Since $\X(g)=0$ the vector $V_g$ writes $\ds V_g=\sum_{i=1}^{n-1}a_iY_i(g)$ in an unique way, and consequently $\ds \norm{V_g}_g=1=\norm{V}_e$.

Let now $g\in \Z\setminus \Z_\X$. At such a point $\ds \X(g)=\sum_{i=1}^{n-1}b_iY_i(g)\neq 0$. Let $\ds b=\sqrt{\sum_{i=1}^{n-1}b_i^2}>0$ and let us choose the particular $\ds V=\sum_{i=1}^{n-1}a_iY_i$ where $a_i=b^{-1}b_i$. Then $V_g=TL_g.V$ writes $\ds V_g=\sum_{i=1}^{n-1}a_iY_i(g)$ but also as
$$
V_g=v\X_g+\sum_{i=1}^{n-1}u_iY_i(g) \quad \mbox{ with } \quad vb_i+u_i=a_i, \qquad i=1,\dots,n-1.
$$
The square of the norm of $V_g$ is $\ds \inf\{v^2+\sum_{i=1}^{n-1}u_i^2;\ vb_i+u_i=a_i \mbox{ for } i=1,\dots,n-1\}$. A straightforward computation shows that the minimum is attained for $\ds v=\frac{b}{1+b^2}$ and consequently that $\ds \norm{V_g}_g^2=\frac{1}{1+b^2}<1$ though $\norm{V}_e=1$.

This shows that the function $g\mapsto\norm{TL_g.V}$ is not continuous at $g$.

Therefore the sub-Riemannian metric is continuous at points of $\Z_\X$ (and at Riemannian points of course) but discontinuous at the points of $\Z\setminus \Z_\X$. This implies that the image of the set $\Z\setminus \Z_\X$ by the isometry $\Phi$ should be $\Z'\setminus \Z'_{\X'}$, which is the desired conclusion.

\hfill $\blacksquare$

\vskip 0.2 cm

\newtheorem{LeftTranslation}[Prop1]{Proposition}
\begin{LeftTranslation} \label{LeftTranslation}
Let $g\in G$. The left translation $L_g$ is an isometry of \sig if and only if $g$ belongs to the set $\Z_\X$ of fixed points of $\X$.
\end{LeftTranslation}

\demo

According to Proposition \ref{Prop1} the condition is necessary.

Conversely let $g\in\Z_\X$. Since $Y_1,\dots,Y_{n-1}$ are left-invariant we have $(L_g)_*Y_i=Y_i$ for $i=1,\dots,n-1$. On the other hand $(L_g)_*\X=\X+Z$ for some right-invariant vector field $Z$. Indeed
$$
\forall Y\in\g \qquad [Y,\X]=(L_g)_*[Y,\X]=[(L_g)_*Y,(L_g)_*\X]=[Y,(L_g)_*\X].
$$
This proves that $(L_g)_*\X$ is an affine vector field whose adjoint action on $\g$ is the same as the one of $\X$. According to \cite{Jou09} it is equal to $\X$ up to a right-invariant vector field $Z$.

Moreover $(\X+Z)(g)=TL_g.\X_e=0$. But $\X_g=0$ and necessarily $Z=0$. This proves that $(L_g)_*\X=\X$, and together with $(L_g)_*Y_i=Y_i$ for $i=1,\dots,n-1$, that $L_g$ is an isometry.

\hfill $\blacksquare$

\vskip 0.2 cm

\newtheorem{Th1}{Theorem}
\begin{Th1} \label{Th1}
The ARSs \sig and \sigp are isometric if and only if there exists an isometry $\Phi$ from \sig onto \sigp such that $\Phi(e)=e$.
\end{Th1}

\demo

Let $\Psi$ be an isometry from \sig onto \sigp. Firstly $\Psi(e)\in \Z'_{\X'}$ according to Proposition \ref{Prop1}. According to Proposition \ref{LeftTranslation} the left translation $L_{\Psi(e)^{-1}}$ is an isometry of \sigp. Then $\Phi=L_{\Psi(e)^{-1}}\circ \Psi$ suits.

\hfill $\blacksquare$

\subsection{Preservation of the distribution}

We recall the following notations: if $Y\in T_eG$ and $g\in G$ then $Y_g$ stands for $TL_g.Y\in T_gG$. 

The key point of this subsection is the following lemma that characterizes the elements of $\Delta_g$ (when $g$ does not belong to the singular locus): they are the vectors whose norm is invariant under (small) left translations.

\newtheorem{CharDelta}{Lemma}
\begin{CharDelta}\label{CharDelta}
Let $\Sigma$ be an ARS on the Lie group $G$. There exists an open and dense subset of $G\setminus \Z$, hence of $G$, of points $g$ that satisfy:
$$
\Delta_g=\{Y_g;\ \norm{Y_h}_h=\norm{Y_g}_g \ \ \mbox{in a neighborhood of }g \}.
$$
\end{CharDelta}

\demo Let $g\in G\setminus \Z$ and $Y\in T_eG$. If $Y=\sum_{i=1}^{n-1}a_iY_i\in\Delta$, then $Y_h$ writes uniquely as $Y_h=\sum_{i=1}^{n-1}a_iY_i(h)$ and $\norm{Y_h}_h$ is constant in a neighborhood of $g$.

We have to show that this is no longer true if $Y\notin\Delta$.

Since the singular locus is an analytic set, we can find $X\in T_eG$ such that $\exp(tX)\notin \Z$ for $t$ in some open interval $(0,\tau)$. Recall that $F$ stands for the mapping from $G$ to $T_eG$ defined by $F(g)=TL_{g^{-1}}.\X_g$. The previous condition is equivalent to $F(\exp(tX))\notin \Delta$ for $t\in (0,\tau)$, or better to:
\begin{equation} \label{toto}
\langle\omega,F(\exp(tX))\rangle\neq 0 \qquad \mbox{ for }\ \ t\in (0,\tau)
\end{equation}
where $\omega$ stands for a nonvanishing one-form on $\g$ orthogonal to $\Delta$.
But Condition (\ref{toto}) together with  $\langle\omega,F(e)\rangle=0$ imply
$$
\dt \langle\omega,F(\exp(tX))\rangle\neq 0
$$
on $(0,\tau)$ (except may be at isolated points). Consequently the set
$$
\Omega=\{g\in G;\ \ g\notin\Z \mbox{ and } \dto \langle\omega,F(g\exp(tX))\rangle\neq 0 \}
$$
is not empty. This set being defined by analytic conditions is open and dense in $G$.

Let $\g\in\Omega$, let us consider the path $\gamma(t)=F(g\exp(tX))$ in $T_eG$ and let $\ds Y=\gamma(0)=F(g)$. Since $g\notin\Z$ the vector $F(g)$ does not belong to $\Delta_e$, the $n$-uple $(F(g),Y_1,\dots,Y_{n-1})$ is a basis of $T_eG$, and the derivative $\gamma'(0)$ can be written as $\ds \gamma'(0)=vF(g)+\sum_{i=1}^{n-1}u_iY_i=vY+\sum_{i=1}^{n-1}u_iY_i$. Then we have
$$
\begin{array}{rll}
\gamma(t) & =Y+t\gamma'(0)+O(t^2)) & \mbox{hence}\\
Y & =\gamma(t)-t\gamma'(0)+O(t^2) & \\
  & =\gamma(t)-t(v\gamma(0)+\sum_{i=1}^{n-1}u_iY_i)+O(t^2) & \mbox{and}\\
(1+tv)Y & =  \gamma(t)-t\sum_{i=1}^{n-1}u_iY_i+O(t^2).
\end{array}
$$
Let us denote by $\norm{.}_t$ the Euclidean norm in $T_eG$ defined by the moving orthonormal frame $(\gamma(t), Y_1,\dots,Y_{n-1})$.
We have
$$
\ds (1+tv)^2\norm{Y}_t^2=1+O(t^2) \quad \mbox{ and } \quad \norm{Y}_t^2=\frac{1}{(1+tv)^2}+O(t^2).
$$
From that we obtain
$$
\dto \norm{Y}_t^2=-2v.
$$
Since $g\in\Omega$ the vector $\gamma'(0)$ does not belong to $\Delta$ and $v\neq 0$. Moreover we have obviously $\norm{Y_{\gamma(t)}} _{\gamma(t)}=\norm{Y}_t$, which proves that in any neighborhood of $g$, $\norm{Y_h}_h$ is not constant.

The same happens for any $Z$ that does not belong to $\Delta$ because it is the sum of an element of $\Delta$ and of $\alpha Y$ for some nonzero real number $\alpha$.

\hfill $\blacksquare$
\newtheorem{preservation}[Th1]{Theorem}
\begin{preservation} \label{preservation}
Let $\Phi$ be an isometry from \sig onto \sigp that preserves the identity. Then:
\begin{enumerate}
\item Its tangent mapping $\Phi_*$ sends $\Delta$ on $\Delta'$, that is, $T_g\Phi.\Delta_g=\Delta_{\Phi(g)}$ for all $g\in G$.
\item Either $\Phi_*\X=\X'$, and $T_e\Phi\circ D=D'\circ T_e\Phi$\\
or $\Phi_*\X=-\X'$, and $T_e\Phi\circ D=-D'\circ T_e\Phi$.
\end{enumerate}
\end{preservation}

\demo According to Lemma \ref{CharDelta} there is an open and dense subset $\Omega$ of points $g$ of $G\setminus\Z$ where the distribution is characterized by:
$$
\Delta_g=\{Y_g;\ \norm{Y_h}_h=\norm{Y_g}_g \ \ \mbox{in a neighborhood of }g \}.
$$
This characterization should be preserved by the isometry $\Phi$ hence $T_g\Phi(\Delta_g)=\Delta'_{\Phi(g)}$. The density of $\Omega$ implies that this equality actually holds at all points.

Let now $g\in G\setminus\Z$. The image by $\Phi$ of the orthonormal frame $(\X(g),Y_1(g),\dots,Y_{n-1}(g))$ is the orthonormal frame $(T_g\Phi.\X(g),T_g\Phi.Y_1(g),\dots,T_g\Phi.Y_{n-1}(g))$. But the fact that $T_g\Phi.\Delta_g=\Delta'_{\Phi(g)}$ implies that $T_g\Phi.\X(g)$ is an unitary vector orthogonal to $\Delta'_{\Phi(g)}$, that is,
\begin{equation}
\label{sign} T_g\Phi.\X(g)=\pm \X'(\Phi(g)).
\end{equation}

In this equality the sign is constant on any connected subset of $G$ where $\X$ does not vanish. On the other hand the set $\Z_\X$ of singularities of $\X$ is a subgroup of $G$ whose Lie algebra is $\ker(D)$. The set $G\setminus\Z_\X$ can be disconnected only if the dimension of $\ker(D)$ is $n-1$. In that case let $Y$ be any element of $T_eG$ that does not belong to $\ker(D)$. Then $\exp(tY)$ does not belong to $\Z_\X$, for $t$ in some interval $(0,T)$, with $T>0$. The sign of the equality $T_{\exp(tY)}\Phi.\X(\exp(tY))=\pm \X'(\Phi(\exp(tY)))$ being constant for $t\in(0,T)$, we denote it by $\epsilon(Y) =\pm 1$. We push this equality to $T_eG$ and we derivate it at $t=0$. On the one hand we get (here $\varphi_t$ stands for the flow of $\X$ and we use Formula (\ref{equationfonda})):

$$
\begin{array}{l}
\ds \dto TL_{(\Phi(\exp(tY)))^{-1}}.T_{\Phi(\exp(tY))}.\X(\Phi(\exp(tY)))
\ds = \dto\ \dso \ (\Phi(\exp(tY)))^{-1}\Phi(\varphi_s(\exp(tY)))\\
\\
\quad \ds = \dto\ \dso \ (\Phi(\exp(tY)))^{-1}\Phi(\exp(e^{sD}tY))
\ds = \dso(-T_e\Phi.Y+T_e\Phi e^{sD}Y )\\
\\
\ds \quad =T_e\Phi DY,
\end{array}
$$
and on the other hand:
$$
\dto TL_{(\Phi(\exp(tY)))^{-1}}.\X'_{\Phi(\exp(tY)}=\dto F'(\Phi(\exp(tY))=D'T_e\Phi.Y
$$
according to Formula (\ref{TgF}) recalled in Section \ref{BD}.

Finally we have obtained $T_e\Phi DY=\epsilon(Y)D'T_e\Phi.Y$. By linearity we have also $T_e\Phi D(-Y)=\epsilon(Y)D'T_e\Phi.(-Y)$, and consequently $T_{\exp(-tY)}\Phi.\X(\exp(tY))=\epsilon(Y) \X'(\Phi(\exp(-tY)))$ for $t$ small enough. But $\exp(tY)$ and $\exp(-tY)$ are located in different connected components of $G\setminus \Z_\X$ (if different connected components exist), and this proves that the sign is constant in Formula (\ref{sign}).

\hfill $\blacksquare$

\vskip 0.2cm

\noindent \textbf{Remark}. Since the replacement of $\X'$ by $-\X'$ does not modify the ARS \sigp we can always  assume that $\Phi_*\X=\X'$.


\subsection{The tangent mapping of an isometry}

\newtheorem{TangentMap}[Th1]{Theorem}
\begin{TangentMap} \label{TangentMap}
Let $\Phi_1$ and $\Phi_2$ be two isometries from \sig onto \sigp that preserve the identity.

If $T_e\Phi_1=T_e\Phi_2$ then $\Phi_1=\Phi_2$.
\end{TangentMap}

Before proving the theorem let us recall that the normal Hamiltonian of an ARS is
$$
\mathcal{H}=\frac{1}{2}\langle \lambda_g,\X(g) \rangle^2+\frac{1}{2}\sum_{i=1}^{n-1}\langle \lambda_g,Y_i(g) \rangle^2,
$$
where $\lambda_g\in T_g^*G$. At Riemannian points the geodesics are the projections of the integral curves of the Hamiltonian vector field associated to $\mathcal{H}$ (see \cite{ABB14} for more information).

\vskip 0.2cm

\demo
The diffeomorphism $\Phi_2^{-1}\circ \Phi_1$ being an isometry of \sig, it is sufficient to prove that an isometry of \sig that preserves the identity $e$ and whose tangent mapping at $e$ is the identity of $T_eG$ is itself the identity mapping of $G$.

So let $\Phi$ be such an isometry of \sig.

In order to simplify the notation let $\Psi$ be the lift of $\Phi$ to the cotangent space, that is
$$
\Psi(g,\lambda_g)=(\Phi(g),\lambda_g\circ T_{\Phi(g)}\Phi^{-1}).
$$
It is a classical fact that $\Psi$ is a symplectomorphism. According to Proposition 4.51 of \cite{ABB14} we have
$$
\Psi_*\op {\mathcal{H}}=\op{\mathcal{H}\circ \Psi^{-1}}
$$
where $\mathcal{H}$ stands for the normal sub-Riemannian Hamiltonian and $\op {\mathcal{H}}$ for the associated Hamiltonian vector field.

Let $g$ be a Riemannian point of $G$. A straightforward computation shows that:
\begin{equation}
\label{hamilton}
\mathcal{H}\circ \Psi^{-1}(g,\lambda_g)=\frac{1}{2}\langle \lambda_g,T_{\Phi^{-1}(g)}\Phi.\X(\Phi^{-1}(g)) \rangle^2+\frac{1}{2}\sum_{i=1}^{n-1}
\langle \lambda_g,T_{\Phi^{-1}(g)}\Phi.Y_i(\Phi^{-1}(g)) \rangle^2.
\end{equation}
Since $\Phi$ is an isometry, $(T_{\Phi^{-1}(g)}\Phi.\X(\Phi^{-1}(g)),T_{\Phi^{-1}(g)}\Phi.Y_1(\Phi^{-1}(g)),\dots,T_{\Phi^{-1}(g)}\Phi.Y_{n-1}(\Phi^{-1}(g)))$ is an orthonormal basis of $T_gG$ and (\ref{hamilton}) is half the sum of the squares of the coordinates of $\lambda_g$ in the dual basis. It is a standard fact of linear algebra that this is invariant by isometry. Consequently
$$
\mathcal{H}\circ \Psi^{-1}(g,\lambda_g)=\frac{1}{2}\langle \lambda_g,\X(g) \rangle^2+\frac{1}{2}\sum_{i=1}^{n-1}\langle \lambda_g,Y_i(g) \rangle^2=\mathcal{H}(g,\lambda_g).
$$
The set of Riemannian points being dense in $G$, this implies $\Psi_*\op {\mathcal{H}}=\op{\mathcal{H}\circ \Psi^{-1}}=\op {\mathcal{H}}$.

Let $(g(t),\lambda(t))$ be a normal extremal with initial conditions $g(0)=e$, and $\lambda(0)=\lambda_g\in T^*_eG$. Then $\Psi(g(t),\lambda(t))$ is also a normal extremal with initial conditions $e$ and $(T_e\Phi^{-1})^*\lambda_g$. Since $T_e\Phi=Id$ we get $\Phi(g(t))=g(t)$ for all $t$ for which $g(t)$ is defined. This is enough to prove that $\Phi=Id$.

\hfill $\blacksquare$


\subsection{Nilpotent groups}\label{NilpotentGroups}
\newtheorem{Nilpotent}[Th1]{Theorem}
\begin{Nilpotent} \label{Nilpotent}
If the group $G$ is nilpotent then the isometries of ARSs of $G$ that preserve the identity are automorphisms.
\end{Nilpotent}

\demo
Let $(\Sigma)$ and $(\Sigma')$ be two ARSs on $G$, and $\Phi$ an isometry from $(\Sigma)$ onto $(\Sigma')$ that preserves the identity.

(i) \textit{\underline{$\Delta$ and $\Delta'$ are not subalgebras}.} In that case $\Delta$ and $\Delta'$ together with their left-invariant metrics define left-invariant sub-Riemannian metrics on $G$. Since $T_g\Phi$ is an isometry from $\Delta_g$ onto $\Delta_{\Phi(g)}$ for all $g\in G$, the diffeomorphism $\Phi$ is an isometry of these left-invariant structures. According to the theorem of Kivioja-Le Donne (\cite{KL16}), $\Phi$ is an automorphism.

(ii) \textit{\underline{$\Delta$ and $\Delta'$ are subalgebras}.} In that case $\Delta$ and $\Delta'$ define left-invariant integrable distributions on $G$. If we denote by $H$ (resp. by $H'$) the connected Lie subgroup generated by $\Delta$ (resp. by $\Delta'$) then the leaves of the induced foliations are the cosets $gH$ (resp. the cosets $gH'$).

Let $g\in G$, and let $M=\Phi(gH)$ be the image of the coset $gH$ by $\Phi$. Since $\Phi$ is a diffeomorphism, $M$ is a submanifold of $G$, and for any $h\in H$ the tangent space at the point $\Phi(gh)$ is $T_{gh}\Phi.\Delta_{gh}=\Delta'_{\Phi(gh)}$. Consequently $M$ is the leaf through the point $\Phi(g)$ of the distribution $\Delta'$, hence equal to $\Phi(g)H'$, and $\Phi$ is an isometry from $gH$ onto $\Phi(g)H'$ (because $T_g\Phi$ respects the metric between $\Delta_g$ and $\Delta_{\Phi(g)}$ for all $g\in G$).

But the left translations respect the metric of $\Delta$, and the one of $\Delta'$ as well, so that $L_{\Phi(g)^{-1}}\circ\Phi\circ L_g$ is an isometry from $H$ onto $H'$. These two groups are nilpotent, and according to The Kivioja-Le Donne's Theorem again, the restriction of $\Phi$ to $H$ is an automorphism onto $H'$.

Let $Y\in \Delta$, and $\overline{Y}$ be the associated left-invariant vector field of $H$. Since $L_{\Phi(g)^{-1}}\circ\Phi\circ L_g$ is an isometry from $H$ onto $H'$, we get that $(L_{\Phi(g)^{-1}}\circ\Phi\circ L_g)_*\overline{Y}$ is a left-invariant vector field $\overline{Z}$ of $H'$.

Let $h\in H$. Then $\Phi(gh)$ belongs to $\Phi(g)H'$, and is equal to $\Phi(g)h'$ for $h'=L_{\Phi(g)^{-1}}\circ\Phi\circ L_g(h)$. The equality
$$
T_h(L_{\Phi(g)^{-1}}\circ\Phi\circ L_g).\overline{Y}_h=\overline{Z}_{h'}
$$
implies
$$
\begin{array}{l}
T_{gh}\Phi.Y_{gh}=T_{gh}\Phi.T_eL_{gh}.Y_e=T_{gh}\Phi.T_hL_{g}.\overline{Y}_h\\
\qquad\qquad\qquad =T_{h'}L_{\Phi(g)}.\overline{Z}_{h'}=T_{h'}L_{\Phi(g)}T_eL_{h'}.Z_e=T_eL_{\Phi(g)h'}.Z_e=T_eL_{\Phi(gh)}.Z_e.
\end{array}
$$
Since $g$ and $h$ are arbitrary, this proves that the image by $\Phi_*$ of the \textbf{left-invariant vector field $Y$ of $G$} is a left-invariant vector field of $G$.

This being true for all $Y\in\Delta$, and according to the forthcoming lemma \ref{LeftInvariantII} the isometry $\Phi$ is an automorphism.

\hfill $\blacksquare$

\vskip 0.2cm

\noindent\textbf{ Remark}. In the algebra case, the proof uses only that the groups $H$ and $H'$ are nilpotent, not that the group $G$ itself is nilpotent.

\newtheorem{LeftInvariantI}[CharDelta]{Lemma}
\begin{LeftInvariantI}\label{LeftInvariantI}
Let $G$ be a connected Lie group and $\g$ its Lie algebra, identified with the set of left-invariant vector fields.

Let $\Phi$ be a diffeomorphism of $G$ that verifies $\Phi(e)=e$ and $\Phi_* Y\in \g$ for all $Y\in \g$.
Then $\Phi$ is an automorphism.
\end{LeftInvariantI}

\demo Let $Y\in\g$ be a left-invariant vector field and $Z=\Phi_* Y$. The flow of $Y$ is $(t,g)\longmapsto g\exp(tY)$ and the one of $\Phi_*Y$ is $(t,g)\longmapsto \Phi(\Phi^{-1}(g)\exp(tY))$. On the other hand $Z$ is left-invariant by assumption and we get:
\begin{equation}\label{LiLi}
\forall g\in G\ \ \forall t\in\R \quad \Phi(\Phi^{-1}(g)\exp(tY))=g\exp(tZ).
\end{equation}

Formula (\ref{LiLi}) applied at $g=e$ gives $\Phi(\exp(tY))=\exp(tZ)$ for all $t\in\R$.

Let $V$ be a neighborhood of $e$ in $G$ and $U$ be a neighborhood of $0$ in $\g$ such that $\exp$ be a diffeomorphism from $U$ onto $V$. If $Y\in U$, then $y=\exp(Y)\in V$ and for all $g\in G$ holds:
$$
\Phi(\Phi^{-1}(g)y)=\Phi(\Phi^{-1}(g)\exp(Y))=g\exp(Z)=g\Phi(\exp(tY))=g\Phi(y).
$$ 
In other words we have $\Phi(zy)=\Phi(z)\Phi(y)$ for all $z=\Phi^{-1}(g)\in G$ and all $y\in V$. Since $V$ generates $G$ this equality is true for all $z,y\in G$ and $\Phi$ is an automorphism.

\hfill $\blacksquare$

\newtheorem{LeftInvariantII}[CharDelta]{Lemma}
\begin{LeftInvariantII}\label{LeftInvariantII}
Let $\Delta$ be a subspace of $\g$ and $\X$ be a linear vector field such that the Lie algebra generated by $\Delta$ and $\X$ be equal to $\g\oplus \R \X$.

If $\Phi$ is a diffeomorphism  of $G$ such that $\Phi_*Y\in\g$ for all $Y\in \Delta$ and such that $\Phi_*\X$ is a linear vector field, then $\Phi$ is an automorphism.
\end{LeftInvariantII}

\demo 
Since the Lie bracket of a left-invariant vector field with another left-invariant vector field or a linear field is left-invariant,  the Lie algebra generated by the $\Phi_*Y$, $Y\in \Delta$, and $\Phi_*\X$ is obviously equal to $\g\oplus \R\Phi_*\X$. According to Lemma \ref{LeftInvariantI} the diffeomorphism $\Phi$ is an automorphism.

\hfill $\blacksquare$

\vskip 0.2cm

Theorem \ref{BD} is no longer true in general if the group is not nilpotent. Counter-examples can be easily built with the help of the Milnor example of the rototranslation group (see \cite{Milnor76}) recalled in \cite{KL16}.

\vskip 0.2cm

\noindent \textbf{Counter-Example}. 
The rototranslation group is the universal covering of the group of orientation-preserving isometries of the Euclidean plane. It can be described as $\R^3$ with the product:
$$
\begin{pmatrix}x \\ y\\ z\end{pmatrix}.\begin{pmatrix}x' \\ y'\\ z'\end{pmatrix}=
\begin{pmatrix}\cos z & -\sin z & 0 \\ \sin z & \cos z & 0\\ 0 & 0 & 1\end{pmatrix}\begin{pmatrix}x \\ y\\ z\end{pmatrix}+\begin{pmatrix}x' \\ y'\\ z'\end{pmatrix}
$$
The Euclidean metric is left-invariant for this product. It can be shown that the group of automorphisms of this group that are isometries is one-dimensional though the group of isometries that preserve the identity is 3-dimensional, which implies that not all isometries are affine.

Let us call this group $\mathcal{R}$ and let us define an ARS on $G=\mathcal{R}\times \R^2$ in the following way: the structure of $\mathcal{R}$ is the previous one, and the one of $\R^2$ is the Grushin plane. Any diffeomorphism $\Phi$ that preserves the identity, made of the direct product of an isometry of $\mathcal{R}$ and an isometry of the Grushin plane is an isometry of this ARS. However if the isometry of $\mathcal{R}$ is not an automorphism, $\Phi$ cannot be an automorphism of $G$.

Here $G$ is a solvable, not nilpotent group.

Many different counter-examples can be build, for instance by replacing the Grushin plane by one of the ARSs defined later on the Heisenberg group.


\subsection{Conclusion}

The results of this section show that an isometry that preserves the identity also preserves the left-invariant distribution, the linear field and is characterized by its tangent mapping at the identity. Moreover that last should preserve the derivation.

We are consequently interested in diffeomorphisms $\Phi$ of $G$ that satisfy:
\begin{enumerate}
\item[(i)] $\Phi(e)=e$
\item[(ii)] $\Phi_*\Delta=\Delta'$
\item[(iii)] $\Phi_*\X=\X'$
\item[(iv)] $T_e\Phi\circ D=D'\circ T_e\Phi$
\end{enumerate}
 Since the isometries that preserve the identity are completely determined by their tangent maps at the origin we will first look for invertible linear maps $P$ on $T_eG$ that verifies $P\circ D=D'\circ P$ and $P(\Delta_e)=\Delta'_e$.
 
 If the Lie group $G$ is simply connected and if such a $P$ is an automorphism of $\g$, it is the tangent mapping of an automorphism $\Phi$ of $G$. It is easy to see that this automorphism is an isometry. Indeed it transforms any left-invariant vector field into a left-invariant vector field. Since $\Phi$ satisfies $P(\Delta_e)=\Delta'_e$, it satisfies $\Phi_*\Delta=\Delta'$. On the other hand it transforms any linear vector field into a linear vector field. Since it satisfies $P\circ D=D'\circ P$ we have $\Phi_*\X=\X'$.
 
 If either $P$ is not an automorphism of $\g$ or the Lie group is not simply connected we cannot conclude so easily, and we have to look in each case to the existence of an isometry $\Phi$ such that $T_e\Phi=P$.
 
 \vskip 0.2cm
 
 \noindent \textbf{Classification}.
 
 In the two following sections we classify the ARSs on the 2D affine group and the Heisenberg group. This classification is done in three steps:
 \begin{enumerate}
 \item The ARSs are at first classified by isometry.
 \item It is clear that a global rescaling of the metric does not modify the geometry. We consequently normalize the ARSs by rescaling in a second step.
 \item In order to emphasize the main geometric structures we then accept to modify the left-invariant metric in the left-invariant distribution $\Delta$ (but we keep it Euclidean and left-invariant). This amounts to "forget" the metric in $\Delta$. The typically different geometries are thus exhibited.
 
 \end{enumerate}


\section{Classification of the ARSs on the affine group} \label{Aff}
 
Let $G$ be the connected $2$-dimensional affine group:
$$
G=Aff_+(2)=\left\{\begin{pmatrix}x&y\\0&1\end{pmatrix};\ \ (x,y)\in\R_+^*\times \R\right\}.
$$
Its Lie algebra is solvable\footnote{A Lie algebra $\g$ is solvable if its derived series terminates in the zero Lie algebra: define $\mathcal{D}^1\g=[\g,\g]$ and by induction $\mathcal{D}^{n+1}\g=[\mathcal{D}^n\g,\mathcal{D}^n\g]$, then $\g$ is solvable if $\mathcal{D}^n\g$ vanishes for some integer $n$.} and generated by the left-invariant vector fields: 
$$
gX=\begin{pmatrix}x&0\\0&0\end{pmatrix} \quad \mbox{ and }\quad gY=\begin{pmatrix}0&x\\0&0\end{pmatrix}\quad \mbox{ where} \quad g=\begin{pmatrix}x&y\\0&1\end{pmatrix}.
$$
In natural coordinates they write $\ds X(x,y)= x\frac{\partial}{\partial x}$ and 
 $\ds Y(x,y)= x\frac{\partial}{\partial y}.$ They verify\\ $[X,Y]=XY-YX=Y$.

In the basis $(X,Y),$ all derivations $D$ of the Lie algebra $\mathfrak{aff}(2)$ have the form $$D=\begin{pmatrix}0&0\\a&b\end{pmatrix},\text{ where } a,b\in\R .$$ 
 
The linear vector field $\X$ associated to such a derivation is $$\X(g)=\begin{pmatrix}0 &a(x-1)+by\\0 &0\end{pmatrix}.$$ In natural coordinates it writes $\ds \X(x,y)=(a(x-1)+by)\frac{\partial}{\partial y}$.
For more details, see \cite{DJ14}.

An ARS on $Aff_+(2)$ is defined by a left-invariant vector field $B=\alpha X+ \beta Y$ and a derivation $D$ such that $B$ and $DB$ are linearly independent, in order to satisfy the rank condition.
In natural coordinates, the ARS is described as the system 
$$\left\lbrace\begin{array}{cl}
\dot{x}= & u\alpha x \\
\dot{y}= & v(a(x-1)+b y)+u\beta x  
 \end{array} \right. $$

\newtheorem{IsometryAff}[Prop1]{Proposition}
\begin{IsometryAff}\label{IsometryAff}
Let $\Sigma=(\mathcal{X},B),$ $\Sigma'=( \mathcal{X}', B' )$ be two ARSs on $Aff_{+}(2).$ 
If $\Phi$ is an isometry between $\Sigma$ and $\tilde{\Sigma}$ and $\Phi(e)=e$ then $\Phi$ is an automorphism.
\end{IsometryAff}

\demo Note that Theorem \ref{preservation} ensures that $\Phi_{*}(\Delta)=\Delta', \ \Phi_{*}(\mathcal{X})=\epsilon \mathcal{X}',$ where $\epsilon=\pm 1.$

Since $\Phi_{*}$ preserves the Euclidean metric in $\Delta,$ which is  1-dimensional, we have 
$$\Phi_{*}(B)=\epsilon' B', \quad \mbox{with }\ \epsilon'=\pm 1.$$

In other words $\Phi_*Y$ is a left-invariant vector field for all $Y\in\Delta$. Together with $\Phi(e)=e,$ and according to Lemma \ref{LeftInvariantII}, we obtain that $\Phi$ is an automorphism.

\hfill $\blacksquare$

We are now in position to state a complete classification of ARSs on $Aff_{+}(2),$ and to characterize them by geometric invariants. 

\subsection{Classification by isometries}

The almost-Riemannian structures classification by isometries is simplified by two facts. First, Proposition \ref{IsometryAff} states that isometries fixing the identity are Lie group automorphism, and second $Aff_{+}(2)$ is simply connected. Thus, instead of searching for Lie groups automorphism, is enough to search for Lie algebra automorphisms.

Let $P$ be an automorphism of the Lie algebra $\mathfrak{aff}(2).$ It is easy to see that $P$ has form
$$P=\left(\begin{array}{cc}
 1 & 0 \\
 c & d
 \end{array} \right), \text{\ for\ } c\in\mathbb{R} \text{ and } d\in\mathbb{R}^{*}$$ and is associated to the Lie group automorphism $\Phi(x,y)=(x,c(x-1)+dy),$ via the equality $P=T_{e}\Phi.$

\newtheorem{ClasAff2}[Prop1]{Proposition}
\begin{ClasAff2}\label{ClasAff2}
Any almost-Riemannian structure on the Lie group $Aff_{+}(2)$ is isometric to one and only one of the structures defined by
$$D=\left(\begin{array}{cc}
 0 & 0 \\
 1 & b
 \end{array} \right)\ \ \ \text{ and }\ \ \ \  B = \alpha X,$$
where $\alpha > 0$ and $b\geq 0.$
\end{ClasAff2}

\demo First note that the change of the vector fields $\mathcal{X}$ to $-\X$ or $B$ to $-B$ does not modify the metric. 

The ARS is defined by $B=\alpha X + \beta Y$ and $D=\begin{pmatrix}0&0\\a&b\end{pmatrix},$ and  we have 
$$DB=D(\alpha X+\beta Y)=\alpha DX+\beta DY=(\alpha a+\beta b)Y.$$
Thus, the 2x2-matrix whose first (resp. second) column contains the coefficients of
$B$ (resp. $DB$) in the basis $\{X, Y\}$ is given by
$$(B\ DB) = \left(\begin{array}{cc}
 \alpha & 0 \\
 \beta & (a\alpha+b\beta)
 \end{array} \right).$$

The rank condition is satisfied if and only if $\alpha(a\alpha+b\beta)\neq 0.$ 

Since $B$ and $DB$ are linearly independent, we can define an isomorphism $P$ of the vector space $\g$ by $$P(B) = \alpha X\ \ \ \ \ \ \text{and}\ \ \ \ \ \ P(DB) = \alpha Y.$$
This isomorphism turns out to be a Lie algebra automorphism on $\g.$ Indeed, the equality
$$[B,DB] = [\alpha X + \beta Y, (a\alpha + b\beta )Y] = \alpha (a\alpha + b \beta )Y = \alpha DB,$$
implies
$$P[B,DB] = P(\alpha DB) = \alpha^{2}Y = \alpha^{2}[X,Y] = [\alpha X, \alpha Y] = [PB,P(DB)].$$

Using Theorem \ref{preservation} the derivation $D$ is transformed to the derivation $PDP^{-1},$ which is characterized by :
$$PDP^{-1}X = Y\ \ \ \ \ \text{and}\ \ \ \ \ PDP^{-1}Y = bY.$$

Therefore, any ARS is isometric to  
$$PDP^{-1}=\left(\begin{array}{cc}
 0 & 0 \\
 1 & b
 \end{array} \right)\ \ \ \text{ and }\ \ \ \  PB = \alpha X$$
with $\alpha > 0$ and $b\geq 0.$ It is easy to check that two different such ARSs are not isometric.

\hfill $\blacksquare$

\

\textbf{Remark:} The singular locus $\Z=\Z_{\X}$ is a Lie subgroup of $Aff_{+}(2).$

\newtheorem{ClasAff3}[Prop1]{Proposition}
\begin{ClasAff3}\label{ClasAff3}
The group of isometries of an ARS on $Aff_{+}(2)$ is the group of left translations by elements of $\Z_{\X}.$
\end{ClasAff3}
\demo
Let $\Phi$ be an isometry. Suppose that $\Phi$ fixes the identity. Proposition \ref{IsometryAff} ensures that any isometry $\Phi$ on $Aff_{+}(2)$ is an automorphism. Since $\Phi_{*}(\Delta)=\Delta,$ we obtain $T_{e}\Phi=\left(\begin{array}{cc}
 1 & 0 \\
 0 & d
 \end{array} \right)$ and using $\Phi_{*}(\X)=\pm \X,$ we conclude $PDP^{-1}=D$ if and only if $d=1,$ i.e. if and only if $\Phi$ is the identity map. 
 
Thus, the isometry group of an ARS is the group of left translations by elements of $g\in \Z_{\X}.$

\hfill $\blacksquare$

\subsection{Global rescaling}

\ 

We do not change the geometry of the ARS, if we multiply all the vector fields by a common positive constant $\lambda.$ 
Choosing $\lambda$ as $\frac{1}{\alpha},$ we rescale $\alpha X$ to $ X$ and the derivation $D$ to 
$$\lambda D=\left(\begin{array}{cc}
 0 & 0 \\
\frac{1}{\alpha} & \frac{b}{\alpha} 
 \end{array} \right). 
\text{ The ARS is isometric to }\left(\begin{array}{cc}
 0 & 0 \\
1 & \frac{b}{\alpha} 
 \end{array} \right),$$
for some $\frac{b}{\alpha}\geq 0.$
Therefore, any ARS is up to a rescaling isometric to an ARS defined by 
$$B= X \text{ and  } D=\left(\begin{array}{cc}
 0 & 0 \\
1 & b 
 \end{array} \right), \quad\mbox{with}\ \  b\geq 0.$$ 
 
Equivalently, it is defined by 
$$B= X \text{ and the singular locus } \Z_{\X}=\{(x,y):(x-1)+by=0, \text{ for some } b\geq0 \}.$$

Moreover, the singular locus is normal if and only if $b=0,$ that is, $\Z_{\X}=\{(1,y), y\in\R \}.$ 
For the other cases, that is, if $b>0$  the singular locus $\Z_{\X}=\{(x,\frac{(x-1)}{b}), x\in\R \} $ is not normal.

\subsection{Deformation of the Euclidean metric in $\Delta$}

\

Deforming the left-invariant metric, we can assume $B=b X,$ (if $b>0$) and the derivation $$ D=\left(\begin{array}{cc}
 0 & 0 \\
b & b 
 \end{array} \right), 
\text{ which are isometric up to a rescale to } B=b X \text{ and }\left(\begin{array}{cc}
 0 & 0 \\
1 & 1 
 \end{array} \right).$$

Therefore, any almost-Riemannian structure on the Lie group $Aff_{+}(2)$ is related by $B=X$ and the derivation $$ D=\left(\begin{array}{cc}
 0 & 0 \\
1 & b 
 \end{array} \right), \text{ where }b=0,1.$$  

Moreover, the singular locus $\Z$ is a normal Lie subgroup if and only if $b=0$ and not normal otherwise.

\subsubsection{Conclusion} 

There two main models, the ones obtained above after rescaling and modification of the left-invariant metric, are completely characterized by the singular locus.

The case where the singular locus is a normal subgroup, that is the the case $b=0$, has been completely studied in \cite{AJ16}.

The other one remains to be analyzed.


\section{Classification of the ARSs on the Heisenberg group}\label{Heisenberg}

Let $G$ be the 3-dimensional Heisenberg group: 
$$ G=\left\{\begin{pmatrix}1&x&z \\0&1&y \\0&0&1 \end{pmatrix};\ \ (x,y,z)\in \R^{3} \right\}.
$$

Its Lie algebra $\mathfrak{g}$ is nilpotent\footnote{A Lie algebra $\g$ is nilpotent if its central series terminates in the zero Lie algebra: define $\mathcal{C}^1\g=[\g,\g]$ and by induction $\mathcal{C}^{n+1}\g=[\mathcal{C}^n\g, \g]$, then $\g$ is nilpotent if $\mathcal{C}^n\g$ vanishes for some integer $n$.} and generated by the left-invariant vector fields:
$$
gX=\begin{pmatrix}0&x&0 \\0&0&0 \\ 0&0&0 \end{pmatrix} \mbox{, }\quad gY=\begin{pmatrix}0&0&x \\0&0&1 \\ 0&0&0 \end{pmatrix} \mbox{ and }\quad gZ=\begin{pmatrix} 0&0&z \\0&0&0 \\ 0&0&0\end{pmatrix} \mbox{ where} \quad g=\begin{pmatrix}1&x&z \\0&1&y \\ 0&0&1\end{pmatrix}.
$$

They verify $[X,Y]=Z$ and the other brackets vanish. In the basis $(X,Y,Z),$ all derivations $D$ of $\mathfrak{g}$ have the form 
\begin{equation}
D =\left(\begin{array}{ccc}
a & b & 0 \\
c & d & 0 \\
e & f & a+d
 \end{array} \right), \text{ where } a,b,c,d,e,f\in\R. \nonumber
\end{equation}

The linear vector field $\X$ associated to such a derivation is 
\begin{equation}
\X(x,y,z)=(ax+by)\frac{\partial}{\partial x}+(cx+dy)\frac{\partial}{\partial y}+(ex+fy+(a+d)z+\frac{1}{2}cx^{2} +\frac{1}{2}by^{2})\frac{\partial}{\partial z}. \nonumber
\end{equation} 
For more details, see \cite{DJ14}.

An ARS on $G$ is defined by an orthonormal frame $\{B_{1},B_{2},\X\},$ where $B_{1},B_{2}$ are left-invariant vector fields and $\X$ is a linear one with associated derivation $D.$ 

To classify the ARSs by isometries, note that Theorem \ref{Nilpotent} states that all isometries that fix the identity are Lie group automorphisms. Since the Heisenberg group is simply connected, is enough to work with Lie algebra automorphisms. 

As we will see, there are two very different cases according to whether $\Delta$  is a subalgebra or not. 
\subsection{$\Delta$ is a subalgebra}

\subsubsection{Classification by isometries}

An automorphism of the Lie algebra $\mathfrak{g},$ in the basis $(X,Y,Z)$, has the form: 
\begin{equation}
P =\left(\begin{array}{ccc}
\alpha & \beta & 0 \\
\gamma & \delta & 0 \\
\epsilon & \zeta & det(A) 
 \end{array} \right), \text{
for A equal to } 
 \left(\begin{array}{cc}
\alpha & \beta \\
\gamma & \delta
\end{array} \right), \nonumber\mbox{ and } \det(A)\neq 0.
\end{equation}  

Note that $Z$ belongs to $\Delta.$  
Indeed, suppose $B_{1}= a_{1}X + b_{1}Y + c_{1}Z$ and $B_{2}=a_{2}X+  b_{2}Y + c_{2}Z,$ the first two coordinates are linearly dependent, because $\Delta$ is a subalgebra. Thus there exists $\lambda\in\R$ such that $(a_{1},b_{1})= \lambda(a_{2},b_{2}).$ If $\lambda=0$ than $B_{2}=c_{2} Z,$ otherwise take the linear combination $ X- \lambda Y=(c_{1}- \lambda c_{2})Z .$  

So, we can choose the orthonormal frame of $\Delta$ as $\{B_{1},\eta Z\},$ for some $\eta \in \R^{*}.$ 

\newtheorem{ClassubH1}[Prop1]{Proposition}
\begin{ClassubH1}\label{ClasH1}Any almost-Riemannian structure on $G,$ whose  distribution $\Delta$ is a subalgebra, is isometric to an almost-Riemannian structure whose orthonormal frame is $\{X,Z,\X\}$ and the derivation $D$ has the following form: $$D=\left(\begin{array}{ccc}
0 & b & 0 \\
c & d & 0 \\
0 & f & d 
\end{array} \right),$$
for some $c>0$ and $d,f\geq 0.$
Moreover, two different ARSs of this form are not isometric.
\end{ClassubH1}

\demo
Let $P$ be the isomorphism defined by $P(B_{1})=X,$ $P(DB_{1})=\frac{a_{1}}{\eta}Y$ and $P(\eta Z)=Z.$ 
The isomorphism $P$ is an automorphism because
$$P[B_{1},Z]=0=[X,\eta Z]=[P(B_{1}),P(Z)].$$
$$P[B_{1},Y]=P(a_{1}Z)=\frac{a_{1}}{\eta} Z=[X,\frac{a_{1}}{\eta}Y]=[P(B_{1}),P(DB_1)].$$

The derivation $\widetilde{D}= PDP^{-1}$ associated to
$X$ verifies $\widetilde{D}X = PDP^{-1}X = PDB = \frac{a_{1}}{\eta}Y.$
Hence, in relation to the basis $(X,Y,Z)$ it has the form $$\widetilde{D}=\left(\begin{array}{ccc}
0 & b & 0 \\
c & d & 0 \\
0 & f & d 
\end{array} \right)\ \text{ where } c=\frac{a_{1}}{\eta}\neq 0.  $$

The changes $\X$ to $-\X$ and $B$ to $-B$ allow to assume $d\geq0$. 
Take 
$$P_{\epsilon\epsilon'}=\left(\begin{array}{ccc}
\epsilon & 0 & 0 \\
0 & \epsilon^{-1}\epsilon' & 0 \\
0 & 0 & \epsilon' 
\end{array} \right),\text{ where }\epsilon,\epsilon'=\pm 1.$$ 
We can apply $P_{\epsilon\epsilon'}$ to change $b$ to $\epsilon\epsilon'b,$ $c$ to $\epsilon\epsilon' c$ and $f$ to $\epsilon f.$ Thus, we can assume $c >0$ and $d,f\geq 0$. It is easy to see that two ARSs related to different derivations of this form are not isometric.  

\hfill $\blacksquare$

\newtheorem{Isometrygroup}[Prop1]{Proposition}
\begin{Isometrygroup}\label{Isometrygroup}
The group of isometries of an ARS in the form of Proposition \ref{ClasH1} is: 

If $d\neq0$ $f\neq0$ the ARS group of isometries is composed only of left translations by elements of $\Z_{\X}.$ 

If $d\neq0$ $f=0$ the ARS group of isometries is composed of left translations by elements of $\Z_{\X}$ and infinitesimal isometries as $P_{\epsilon\epsilon'}$ for $ \epsilon=\pm 1$ and $\epsilon'=1.$

If $d=0$ $f\neq0$ the ARS group of isometries is composed of left translations by elements of $\Z_{\X}$ and infinitesimal isometries as $P_{\epsilon\epsilon'}$ for $ \epsilon=\epsilon'=\pm1.$

If $d=f=0$ the ARS group of isometries is composed of left translations by elements of $\Z_{\X}$ and infinitesimal isometries as $P_{\epsilon\epsilon'}$ for $ \epsilon,\epsilon'=\pm 1.$

Generically $b,d$ and $f$ are nonzero and the set $\Z_{\X}$ is reduced to the identity. Thus in the subalgebra case, the group of isometries is generically reduced to the identity.
\end{Isometrygroup}

\demo Let $\Phi $ be an isometry fixing the identity. It is an automorphism. Let $T_{e}\Phi=P_{\epsilon\epsilon'}.$ Note that $P_{\epsilon\epsilon'}(X)=\epsilon X$ and $P_{\epsilon\epsilon'}(Z)=\epsilon' Z$. We obtain that $P_{\epsilon\epsilon'}$ changes $D$ to 
\begin{equation}\label{eq1}
P_{\epsilon\epsilon'}DP^{-1}_{\epsilon\epsilon'}=\left(\begin{array}{ccc}
0 & \epsilon'b & 0 \\
\epsilon'c & d & 0 \\
0 & \epsilon f & d 
\end{array} \right).
\end{equation}
 Since $\Phi_{*}(\X)=\pm\X$ we have $P_{\epsilon\epsilon'} D(P_{\epsilon\epsilon'})^{-1}=\pm D.$ 

Generically $b,c,d$ and $f$ are non zero, so $\epsilon$ and $\epsilon'$ need to be equal to $1.$ Moreover a straightforward computation proves that $\Z_{\X}$ is reduced to $\{e\}.$ As a consequence the group of isometries is generically reduced to the identity.

Considering Equation (\ref{eq1}) it is easy to conclude in the other cases.

\hfill $\blacksquare$

\subsubsection{Global rescaling}

\
We do not change the geometry of the ARS if we multiply all the vector fields by a common positive constant $\lambda.$ 
This global rescaling allows us to use 
$$P =\left(\begin{array}{ccc}
\alpha & 0 & 0 \\
0 & \epsilon & 0 \\
0 & 0 & \alpha \epsilon 
 \end{array} \right) \text{ for $ \alpha\in \R^{*} $ and $\epsilon=\pm 1 ,$
to changes $D$ into } \widetilde{D}=\left(\begin{array}{ccc}
0 & \alpha \epsilon^{-1} b & 0 \\
\alpha^{-1} \epsilon c & d & 0 \\
0 & \alpha f & d 
\end{array} \right).$$
 
Therefore, any ARS is isometric up to a rescale to one and only one ARS defined by the orthonormal frame $\{X,Z,\X\}$ where the associated derivation $D$ is equal to: 

$$ D=\left(\begin{array}{ccc}
0 & b & 0 \\
1 & d & 0 \\
0 & f & d 
\end{array} \right)   
$$

\subsubsection{Deformation of the Euclidean metric in $ \Delta$}
 
In order to obtain a classification with less parameters, we consider as equivalent two ARSs that have the same linear fields, same distributions but provided with a different left-invariant metrics.  

In other words, we forget the left-invariant metrics.
Let $\Sigma$ be the ARS defined by the orthonormal frame $\{X,Z,\X\},$ where the associated derivation is $D=\left(\begin{array}{ccc}
0 & b & 0 \\ 
1 & d & 0 \\
0 & f & d 
\end{array} \right).  
$ We consider as equivalent to $\Sigma,$ any ARS obtained by rescaling and conjugation by automorphisms that preserve $\Delta,$ that is: 
$$P =\left(\begin{array}{ccc}
\alpha & \beta & 0 \\
0 & \delta & 0 \\
\epsilon & \zeta & \alpha \delta 
 \end{array} \right).$$
Thus, we can conclude that, any almost-Riemannian structure on $G,$ whose left-invariant distribution $\Delta$ is a Lie subalgebra, is related to one and only one ARS whose orthonormal frame is $\{X,Z,\X\},$ and as previously the associated derivation $D$ has the form: 

$$D=\left(\begin{array}{ccc}
0 & b & 0 \\ 
1 & d & 0 \\
0 & f & d 
\end{array} \right),  
$$
but the constants $b$, $d$ and $f$ are normalized as follows:

\vskip 0.4cm

\begin{tabular}{|c|c|c|c|c|c|c|c|}
\hline 
case& d & b & f & singular locus $\Z$ & $\Z_{\mathcal{X}}$  & eigenvalues $l_{1}, l_{2}$ \\  
\hline 
(i) & $1$ & $>\frac{-1}{4}$ & $0$ &  x+y=0   & (0,0,0)   & $l_{1}=\frac{1+\sqrt{1+4b}}{2}$ $l_{2}=\frac{1-\sqrt{1+4b}}{2}$ \\ 
\hline 
(ii) & $1$ & $<\frac{-1}{4}$ & $0$ &  x+y=0 &  (0,0,0)    & $l_{1}=\frac{1+\sqrt{-(1+4b})i}{2}$ $l_{2}=\frac{1-\sqrt{-(1+4b})i}{2}$\\ 
\hline 
(iii) & $1$ & $0$ & $1$ &  x+y=0  & $(-y,y,-y-\frac{1}{2}y^{2})$  & $l_{1}=0$  $l_{2}=1$\\ 
\hline 
(iv)&$1$ & $0$ & $0$ & x+y=0 & $(-y,y,-\frac{1}{2}y^{2})$  & $l_{1}=0$ $l_{2}=1$\\ 
\hline 
(v)&$0$ & $1$ & $0$ & x=0 &  x=y=0 &  $l_{1}=1$ $l_{2}=-1$ \\
\hline 
(vi)& $0$ & $-1$ & $0$ & x=0 & x=y=0 &  $l_{1}=i$ $l_{2}=-i$\\ 
\hline 
(vii) & $0$ & $0$ & $1$ & x=0 & x=y=0 & $l_{1}=l_{2}=0$ \\ 
\hline 
(viii) & $0$ & $0$ & $0$ & x=0 & x=0 & $l_{1}=l_{2}=0$\\ 
\hline 
\end{tabular} 

\vskip 0.4cm

\textbf{Remark:} For this case, we know that $\Z$ is always a Lie subgroup and its Lie algebra $T_{e}\Z= D^{-1}\Delta.$ The Lie subalgebra $\Delta$ is not tangent to the singular locus. 
Therefore, there are no tangency points when $\Delta$ is a subalgebra.

\subsection{$\Delta$ is not a subalgebra}

In this case the distribution $\Delta$ generates the Lie algebra. So, to satisfy the rank condition there are no restriction on the derivation $D.$ On the other hand, Condition (iii) of Definition \ref{DefARS} states that $\X$ and $\Delta$ need be linearly independent in an open and dense set. 
 
\subsubsection{Classification by isometries}

\newtheorem{ClasnosubH1}[Prop1]{Proposition}
\begin{ClasnosubH1}\label{ClasnosubH1}
Any ARS, whose distribution $\Delta$ is not a subalgebra, is isometric to an ARS whose orthonormal frame is $\{X,Y,\X\},$ where the derivation $D$ has the following form:
  $$D=\left(\begin{array}{ccc}
a & b & 0 \\
c & d & 0 \\
0 & f & a+d 
\end{array} \right) \text{with } c,f\geq 0.$$
The conditon that $\X(g)$ does not belong to $\Delta(g)$ everywhere reduces to:
If $b=c=f=0$ then $a+d\neq 0.$
\end{ClasnosubH1} 
 
\demo
If $\Delta$ is not a subalgebra then
$$[B^{1},B^{2}]=(a_{1}b_{2}-a_{2}b_{1})Z\notin \Delta.$$

Let us choose the isomorphism $P$ defined by  $$P(B^{1})=X,\ \ P(B^{2})=Y\ \ \text{and}\ \ PZ=\frac{1}{(a_{1}b_{2}-a_{2}b_{1})}Z.$$ 

Actually, the isomorphism $P$ is a Lie algebra automorphism
$$[P(B^{1}),P(B^{2})]=[X,Y]=Z=P((a_{1}b_{2}-a_{2}b_{1})Z)=P[B^{1},B^{2}].$$

To finish, we can apply a rotation within $\Delta$ in order to simplify the derivation $D$. Consider $$P_{\theta,\epsilon}=\left(\begin{array}{cc}
R_{\theta\epsilon} & 0  \\
0 & det(R_{\theta\epsilon}) 
\end{array} \right)=\left(\begin{array}{ccc}
 \cos(\theta) & -\epsilon \sin(\theta) & 0  \\
 \sin(\theta) & \epsilon \cos(\theta) & 0  \\
 0 & 0 & \epsilon 
\end{array} \right),$$ where $R_{\theta\epsilon}$ is an orthogonal transformation of $\Delta$ and $P_{\epsilon,\epsilon'}$ is defined in previous section.

Case (i) $e,f\neq0,$ $c\in\R.$ Apply $P_{\theta,1}$ for some $\theta$
to vanish $e$ then apply $P_{\epsilon,\epsilon'}$ to transform $f$ to $\epsilon f$ and $c$ to $\epsilon' c.$ So, we can assume $c,f\geq 0.$

Case (ii) $e=0,$ $f\neq 0$ and  $c\in\R.$ Apply $P_{\epsilon,\epsilon'}$ to transform $f$ to $\epsilon f$ and $c$ to $\epsilon' c.$ So, we can assume $c,f\geq 0.$

Case (iii) $e,f=0$ and $c\in\R.$ Apply $P_{\epsilon,\epsilon'}$ to transform $c$ to $\epsilon' c$ so we can assume $c\geq0.$

\hfill $\blacksquare$  

\newtheorem{Isometrygroup2}[Prop1]{Proposition}
\begin{Isometrygroup2}\label{Isometrygroup2}
The group of isometries of such an ARS is generated by the left translations by  elements of $\Z_{\X}$ and some of the automorphisms $P_{\theta,\epsilon}$ defined above. Generically, the group of isometries is reduced to left translations.
\end{Isometrygroup2}

\demo Let $\Phi $ be an isometry fixing the identity. 
Theorem \ref{Nilpotent} ensures that $\Phi$ is an automorphism. Since $\Phi$ preserves the Euclidean metric in $\Delta,$ we have $\Phi_{*}(X)=\alpha X+\epsilon\beta Y,$ $\Phi_{*}(Y)= -\beta X + \epsilon\alpha Y$ and $\Phi_{*}(\X)=\epsilon'\X,$ with 
$\alpha^{2}+\beta^{2}= 1,$ and $\epsilon,\epsilon'=\pm 1.$

Note that the differential $\Phi_{*}$ at the identity has the form of the Lie algebra automorphism $P_{\theta,\epsilon}.$
The derivation $$D=\left(\begin{array}{cc}
A &0 \\
 (e,f) &  tr(A) 
\end{array} \right) \text{ is changed by }P_{\theta\epsilon}\text{ to }P_{\theta\epsilon}DP^{-1}_{\theta\epsilon}=\left(\begin{array}{cc}
R_{\theta\epsilon}AR^{-1}_{\theta\epsilon} &0 \\
\epsilon(e,f)R^{-1}_{\theta\epsilon} &  tr(A) 
\end{array} \right).$$

Note that for $e,f\neq 0,$ $P_{\theta\epsilon'}$ alters $(e,f),$ except when $\theta = 2k\pi$ for $k \in \mathbb{Z}$  and the automorphism  $P_{\epsilon\epsilon'}$ change $b$ to $\epsilon'b,$ $c$ to $\epsilon' c,$ $e$ to $\epsilon \epsilon' c$ and $f$ to $\epsilon f.$ If $(e,f)\neq (0,0)$ to preserve $D$ we need $\epsilon=\epsilon'=1.$ Generically $tr(A),\ e,\ f\neq 0$ and the only solution is $P_{\theta\epsilon'}=I$.

\hfill $\blacksquare$  

\subsubsection{Global rescaling}

We do not change the geometry of the ARS, if we multiply all the vector fields by a common positive constant $\lambda.$ 
This global rescaling allows us to use 
$$P =\left(\begin{array}{ccc}
\alpha & 0 & 0 \\
0 & \alpha & 0 \\
0 & 0 & \alpha^{2} 
 \end{array} \right) \text{ where $\alpha\in \R^{*},$ 
to changes $D$ into } \widetilde{D}=\left(\begin{array}{ccc}
a & b & 0 \\
 c & d & 0 \\
0 & \alpha^{-1} f & a+d 
\end{array} \right).$$
 
Therefore, any ARS is defined up to a rescale by the orthonormal frame $\{X,Z,\X\},$ where the derivation $D$ associated to $\X$ has one of the following form: 
If $f\neq0$ we change it to $1.$ If $f=0$ and $c\neq0$ we change it to 1. If $c,f=0$ and $b\neq0$ we change it to $\pm 1.$ If $b,c,f=0$ and $a+d$ is different from $0,$ so we change it to $1.$ 

\subsubsection{Deformation of the Euclidean metric in $ \Delta$}

In order to get a classification with less parameters, we consider as equivalent ARS that have the same linear fields, same distributions but provided with a different left-invariant metrics.  

In other words, we forget the left-invariant metric on $\Delta.$ Consider the orthonormal frame of the previous section.
  
We consider as equivalent to $\Sigma$ any ARS obtained by rescaling and conjugation by automorphisms that preserve $\Delta,$ that is: 
\begin{equation}
P =\left(\begin{array}{ccc}
\alpha & \beta & 0 \\
\gamma & \delta & 0 \\
0 & 0 & \alpha\delta-\beta\gamma 
 \end{array} \right). \nonumber
\end{equation}

The Lie algebra automorphism $P$ transforms the ARS $(\Sigma)$ to the ARS $(\widetilde{\Sigma})$ defined by the orthonormal frame $\{X,Z,\X\}$ where $\widetilde{\mathcal{X}}$ is associated to the derivation
$$\widetilde{D}=\frac{1}{(\alpha\delta-\beta\gamma)}\left(\begin{array}{ccc}
\delta(\alpha a+\beta c)-\gamma(\alpha b+\beta d) & -\beta(\alpha a+\beta c)+\alpha(\alpha b+\beta d) & 0 \\
\delta(\gamma a+ \delta c)-\gamma(\gamma b+\delta d) & -\beta(\gamma a+\delta c)+\alpha(\gamma b+\delta d) & 0 \\
\delta e -\gamma f & - \beta e + \alpha f & (\alpha\delta-\beta\gamma)(a+d) 
 \end{array} \right).$$

A short computation proves the impossibility to vanish $e$ and $f$ simultaneously. We use the diagonalization process to classify by eigenvalues. By doing so, we may alter the entries $e$ and $f$ again. Thus, first we need to put the derivation $D$ in a block form and then apply the previous computations to obtain the following forms:
 
\begin{equation}
(1)\left(\begin{array}{ccc}
l_{1} & 0 & 0 \\
0 & l_{2} & 0 \\
e & f & l_{1}+l_{2}
 \end{array} \right), 
(2)\left(\begin{array}{ccc}
l_{1} & 1 & 0 \\
0 & l_{1} & 0 \\
e & f & 2l_{1}
 \end{array} \right),
(3)\left(\begin{array}{ccc}
a & b & 0 \\
-b & a & 0 \\
e & f & 2a
 \end{array} \right),
\end{equation}  
where $l_1,l_2,a,b,e,f$ are real numbers, $b\neq 0$ and $a^2+b^2>0.$  
 
For $\Delta=\{X,Y\}$ and $D=\left(\begin{array}{ccc}
a & b & 0 \\
c & d & 0 \\
e & f & a+d
 \end{array} \right)$ the singular locus is described by the following formula
$$\Z=\{ex+fy+(a+d)z-\frac{1}{2}cx^{2}+\frac{1}{2}by^{2}-dxy=0 \}.$$ 
 
\textbf{The case (1)} have four distinct subcases with distinct characteristics:
 
\begin{tabular}{|c|c|c|c|c|c|}
 \hline 
case & $l_1$ & $l_2$ & $l_1 +l_2$ & singular locus $\Z$ &  $\Z_{\mathcal{X}}$  \\ 
 \hline 
$(i)$& $\neq 0$ & $\neq 0$ & $\neq 0$ & $(l_1+ l_2)z=l_{2}xy-ex-fy$ & $x=y=z=0$  \\ 
 \hline 
$(ii)$   & $\neq 0$ & $-l_1$ & 0 & $0=l_{1}xy+ex+fy$ & $x=y=0$  \\ 
 \hline 
 $(iii)$ & $\neq 0$ & $0$ & $l_1$ & $l_1 z=-ex-fy$  & $(0,y,\frac{f}{l_1} y)$ \\ 
 \hline 
$(iv)$  &$0$ & $0$ & 0 & $0=ex+fy$ & $ 0=ex+fy$  \\ 
 \hline 
  \end{tabular} 

\

The case (2) has two different behaviours:
 
 \begin{tabular}{|c|c|c|c|}
 \hline 
 cases & $l_1$  & singular locus $\Z$ & $\Z_{\mathcal{X}}$ \\ 
 \hline 
 $(i)$ & $\neq 0$ & $ex+fy+2l_{1}z+\frac{1}{2}y^{2}-l_{1}xy=0$ & (0,0,0) \\ 
  \hline 
 $(ii)$ & $0$ & $ex+fy+\frac{1}{2}y^{2}=0$ & $(0,0,z)$ \\ 
 \hline 
  \end{tabular}

\

The case $(3)$ has two different behaviours and a restriction $a^{2}+b^{2}>0$:

\begin{tabular}{|c|c|c|c|}
 \hline 
 cases & a & singular locus $\Z$ &  $\Z_{\mathcal{X}}$ \\ 
 \hline 
 $(i)$ & $\neq 0$ & $ex+fy+2az+\frac{1}{2}bx^{2}+\frac{1}{2}by^{2}-axy=0$ & $(0,0,0)$ \\ 
 \hline 
 $(ii)$ & $0$ & $ex+fy+\frac{1}{2}bx^{2}+\frac{1}{2}by^{2}=0$ & $(0,0,z)$ \\ 
  \hline 
\end{tabular}
  
\ 

\textit{Tangency points}.
It is clear that an isometry sends tangency points to tangency points and that to find them it is necessary to describe the singular locus. 

Notice in case $(1.i),$ the singular locus $\Z$ is $g^{-1}(0),$ where the function $g$ is\begin{equation}
g(x,y,z)=-(l_1+ l_2)z+l_{2}xy-ex-fy. \nonumber
\end{equation} 
The differential $dg$ is surjective at all points of the locus that is, $0$ is a regular value of the function $g.$ So, the inverse image $g^{-1}(0)$ is a connected two dimensional manifold. 
This case has one tangency point $(x,y,z)=(\frac{-f}{l_{1}},\frac{e}{l_{2}},\frac{-ef}{l_{2}(l_{1}+l_{2})}).$
 
\textbf{In case $(1.ii)$}, the function $g$ is \begin{equation}
g(x,y,z)=l_{1}xy-ex-fy.
\end{equation} 
The differential $dg$ is not surjective at all points of the singular locus when $f$ or $e$ is equal to zero. 

If $e,f\neq 0,$ the singular locus $\Z$ is a nonconnected two dimensional manifold.
If $f=0,$ the singular locus is the union two planes intersecting in the line 
$(0,\frac{e}{l_{1}},z).$
If $e=0,$ the singular locus is the union two planes intersecting in the line $(\frac{f}{l_{1}},0,z).$ 
If $e,f=0$ the orthonormal frame $\{X,Y,\X\}$ does not describe an ARS.

\textit{Tangency points}. This case has three different behaviours. If $e$ and $f$ are both nonzero, there are no tangency points, if $e=0$ tangency points form a line $(\frac{f}{l_{1}},0,z),$ if $f=0$ tangency points form a line $(0,\frac{e}{l_{1}},z).$ In these cases the respective lines are the intersection of the two planes.   
  
\textbf{In case $(1.iii)$} the function $g$ is \begin{equation}
g(x,y,z)=l_1z+ex+fy.
\end{equation}
Note that the function $g$ describes a plane. So, the singular locus is a connected two dimensional manifold. 

\textit{Tangency points}. This case has two possible behaviours. If $e$ is zero, the tangency points form a 1-dimensional manifold described as $(\frac{-f}{l_{1}},y,\frac{-fy}{l_{1}}),\ y\in \mathbb{R}.$  If $e$ is not zero, there are no tangency points.

\textbf{In case $(1.iv)$} the function $g$ is \begin{equation}
g(x,y,z)=-ex-fy.
\end{equation}
Note that the function $g$ describes a plane. So, the singular locus is a connected two dimensional manifold. If $e$ and $ f$ are both zero then $D$ is the zero derivation, which do not describe an ARS.
There are no tangency points. 

\textbf{In case $(2.i)$} the function $g$ is \begin{equation}
g(x,y,z)=ex+fy+2l_{1}z+\frac{1}{2}y^{2}-l_{1}xy, \nonumber
\end{equation} and $0$ is a regular value. So, the singular locus is a connected two dimensional manifold. 
This ARS has one tangency point given by $(x,y,z)=(-\frac{f}{l_{1}}-\frac{e}{l^{2}_{1}},\frac{e}{l_{1}},-\frac{e^{2}}{4l^{2}_{1}}-\frac{ef}{2l_{1}}).$

\textbf{In case $(2.ii,$} the function $g$ is
\begin{equation}
g(x,y,z)=ex+fy+\frac{1}{2}y^{2}.
\end{equation}
The differential $dg$ is surjective at all points only if $e$ is different from zero.  Under this assumption the singular locus $\Z$ is a connected  two dimensional manifold. 
 
If $e$ is zero, there are many points where the application $dg$ is not surjective. In this case, the set $g^{-1}(0)$ has different behaviours according to the value of $f.$ If $f\neq 0,$ the set $g^{-1}(0)$ is the union two planes perpendicular to the $y-axis,$ one through the point $(0,0,0)$ and the other through the point $(0,-2f,0).$ 
If $f=0,$ the set $g^{-1}(0)$ is a plane perpendicular to the y-axis passing through the point $(0,0,0).$
 
Thus, if $e\neq 0,$ there are no tangency points. If $e=0,$ the tangency points need to satisfy the equation $f+y=0.$ So, for $f\neq 0,$ there are no tangency points. If $f=0$ all the singular locus is equal to the set of tangency points. 

\textbf{In case $(3.i)$} the function $g$ is 
\begin{equation}
g(x,y,z)=ex+fy+2az+\frac{1}{2}bx^{2}+\frac{1}{2}by^{2}-axy.
\end{equation} 
The differential $dg$ is surjective at all points of the singular locus. So, $\Z$ is a connected two dimensional manifold. 
Since the following linear system solution is unique 
\begin{equation}
e+bx-ay=0 \qquad   \text{and} \qquad  f+by+ax=0, \nonumber
\end{equation} 
there is one tangency point. 
 
In case $(3.ii),$ the function $g$ is 
\begin{equation}
g(x,y,z)=ex+fy+\frac{1}{2}bx^{2}+\frac{1}{2}by^{2}.
\end{equation} 
The differential $dg$ is not surjective, but the equality $ex+fy+\frac{1}{2}bx^{2}+\frac{1}{2}by^{2}=0$ represents an elliptic cylinder when $e$ or $f$ is different from $0.$ If $e=f=0,$ it is a line through the identity. 
If $e$ or $f$ is different from $0,$ there are no tangency points.
If $e$ and $f$ both vanish, the tangency points forms the line $(0,0,z)\ z\in\mathbb{R}.$
\newline 

In order to simplify even more the derivation, is possible to apply one or both of the following two kinds of automorphisms:

\begin{equation}
P_{\alpha.\beta}=\left(\begin{array}{ccc}
\alpha & 0 & 0 \\
0 & \beta & 0 \\
0 & 0 & \alpha\beta
 \end{array} \right)\ \text{and}\ P_\theta =\left(\begin{array}{ccc}
\cos(\theta) & -\sin(\theta) & 0 \\
\sin(\theta) & \cos(\theta) & 0 \\
0 & 0 & 1
 \end{array} \right), \nonumber
\end{equation} 
depending on each case, together with a global rescaling. We obtain that
 any almost-Riemannian structure on $G,$ whose left-invariant distribution $\Delta$ is not a Lie subalgebra, is related to only one ARS whose the orthonormal frame is $\{X,Y,\X\}$ where the derivation $D$ has the following form: 

$$ D_{1}=\left(\begin{array}{ccc}
l_{1} & 0 & 0 \\
0 & l_2 & 0 \\
e & f & l_{1}+l_2
 \end{array} \right) \ \ \ 
D_{2}=\left(\begin{array}{ccc}
l_{1} & 1 & 0 \\
0 & l_{1} & 0 \\
e & f & l_{1}
 \end{array} \right)
\ \ \ \
D_{3}=\left(\begin{array}{ccc}
a & -b & 0 \\
b & a & 0 \\
e & f & 2a
 \end{array} \right) 
$$

\begin{tabular}{|c|c|c|c|c|c|c|}
\hline 
case& $l_{1}$ & $l_{2}$ & $(e,f)$ & singular & tangency & number of cc \\  
 &   &   &   & locus $\Z$ & points & of $G-\Z$ \\ 
\hline 
 1.i.1  & $1$ & $\mathbb{R}-\{0,-1\}$ & $(1,1)$ & submanifold   & $(-1 ,\frac{1}{l_{2}},\frac{-1}{l_{2}(1+l_{2})})$ & two \\ 
\hline 
 1.i.2  & $1$ & $\mathbb{R}-\{0,-1\}$ & $(1,0)$ & submanifold & $(0,\frac{1}{l_{2}},0)$ & two \\ 
\hline 
 1.i.3  & $1$ & $\mathbb{R}-\{0,-1\}$ & $(0,0)$ & submanifold & $(0,0,0)$ & two \\ 
\hline 
 1.ii.1  &$1$ & $-1$ & $(1,1)$ & submanifold & no tangency points & three \\ 
\hline 
 1.ii.2  & $1$ & $-1$ & $(1,0)$ & not submanifold & $(0,1,z), z\in \mathbb{R}.$ & four\\ 
\hline 
 1.iii.1  &$1$ & $0$ & $(1,1)$ & submanifold & no tangency points & two\\ 
\hline 
 1.iii.2 &$1$ & $0$ & $(0,1)$ & submanifold & $(-1,y,-y), y\in \mathbb{R}.$ &  two \\ 
\hline 
 1.iii.3 &$1$ & $0$ & $(1,0)$ & submanifold & no tangency points &  two\\ 
\hline 
 1.iii.4 &$1$ & $0$ & $(0,0)$ & submanifold & $x=z=0$ &  two \\ 

\hline 
 1.iv.1 &$0$ & $0$ & $(0,1)$ & Lie subgroup & no tangency points & two \\ 
\hline 
\end{tabular} 

* cc means connected components

\

\begin{tabular}{|c|c|c|c|c|c|}
\hline 
case& $l_{1}$ &  $(e,f)$ & singular & tangency & number of cc \\  
 &   &      & locus $\Z$ & points & of $G-\Z$  \\ 
\hline 
 2.i.1  & $\mathbb{R}^{\ast}$  & $(1,1)$ & submanifold   & $(-2,1,-\frac{3}{4}) $ & two  \\ 
\hline 
 2.i.2  & $1$  & $(0,1)$ & submanifold   & $(-1,0,0) $ & two \\ 
\hline 
 2.i.3  & $1$  & $(1,0)$ & submanifold   & $(-1,1,-\frac{1}{4}) $ & two \\ 
\hline 
 2.i.4  & $1$  & $(0,0)$ & submanifold   & $(0,0,0)$ & two \\ 
\hline 
 2.ii.1  & $\mathbb{R}^{\ast}$  & $(1,1)$ & submanifold   & no tangency points & two \\ 
\hline 
 2.ii.2  & $1$  & $(0,1)$ & submanifold   & no tangency points & three \\ 
\hline 
 2.ii.3  & $1$  & $(1,0)$ & submanifold   & no tangency points & two \\ 
\hline
 2.ii.4  & $1$  & $(0,0)$ & Lie subgroup   & tangency point set is equal to $\Z$ & two \\ 
\hline 
\end{tabular}

\

\begin{tabular}{|c|c|c|c|c|c|c|}
\hline 
case& a & b &  $(e,f)$ & singular & tangency  & number of cc  \\  
 &   &   &      & locus $\Z$ & points & of $G-\Z$\\  
\hline 
 3.i.1  & 1  & $\mathbb{R}^{\ast}$ & $(0,1)$ & submanifold   & one tangency point & two \\ 
\hline 
 3.i.2  & 1  & $\mathbb{R}^{\ast}$ & $(0,0)$ & submanifold   & (0,0,0) & two \\ 
\hline 
 3.ii.1  & 0  & $1$ & $(0,1)$ & submanifold   & no tangency points & two \\ 
\hline 
 3.ii.2  & 0 & $1$ & $(0,0)$ & Lie subgroup   & the line $x=y=0.$ &   one \\ 
\hline 
\end{tabular}







\begin{thebibliography}{999}
\bibitem{ABB14} A. Agrachev, D. Barilari, U. Boscain, \textit{Introduction to Riemannian and Sub-Riemannian Geometry}, http://webusers.imj-prg.fr/~davide.barilari/Notes.php
\bibitem{ABS08} A. Agrachev, U. Boscain, M. Sigalotti, \textit{A Gauss-Bonnet like formula on two-dimensional almost-Riemannian manifolds}, Discrete Contin. Dyn. Syst. 20 (4) (2008) 801-822.
\bibitem{ABCGS10} A. Agrachev, U. Boscain, G. Charlot, R. Ghezzi, and M. Sigalotti, \textit{Two dimensional
almost-Riemannian structures with tangency points}. Ann.
Inst. H. Poincaré. Anal. Non Linéaire 27 (2010), 793-807.
\bibitem{AJ16} V.~Ayala, Ph.~Jouan \textit{Almost-Riemannian Geometry on Lie groups}, SIAM J. Control and Optimization 54 (2016), no.5, 2919-2947.
\bibitem{AS01} L.~San Martin and V.~Ayala \textit{Controllability properties of a class of control systems on Lie Groups}, Nonlinear control in the year 2000, Vol. 1,  83--92, L.N. in Control and I.S., 258, Springer, 2001.
\bibitem{AT99} V.~Ayala and J.~Tirao \textit{Linear control systems on Lie groups and Controllability}, Proceedings of Symposia in Pure Mathematics, Vol 64, AMS, 1999, 47-64.
\bibitem{BCST09} B. Bonnard, J.-B. Caillau, R. Sinclair, M. Tanaka, \textit{Conjugate and cut loci of a two-sphere of revolution with application to optimal control},
Ann. Inst. H. Poincaré Anal. Non Linéaire 26 (4) (2009) 1081-1098.
\bibitem{BCGJ11} B. Bonnard, G. Charlot, R. Ghezzi, G. Janin, \textit{The sphere and the cut locus at a tangency point in two-dimensional almost-Riemannian geometry}, J. Dyn.
Control Syst. 17 (1) (2011) 141-161.
\bibitem{BCGM14} U.~Boscain, G.~Charlot, M.~Gaye, P.~Mason \textit{Local properties of almost-Riemannian structures in
dimension 3}, Discrete and Continuous Dynamical Systems-A, Volume 35, Issue 9 (2015), pp. 4115-4147.
\bibitem{BCGS} U. Boscain, G. Charlot, R. Ghezzi, M. Sigalotti, \textit{Lipschitz classification of almost-Riemannian distances on compact oriented surfaces}, J. Geom. Anal.,January 2013, Volume 23, Issue 1, pp 438-455.
\bibitem{BCG13}U. Boscain, G. Charlot, R. Ghezzi \textit{Normal forms and invariants for 2-dimensional almost-Riemannian
structures}, Differential Geometry and its Applications 31 (2013) 41-62.
\bibitem{DJ14}  M.~Dath, Ph.~Jouan \textit{Controllability of linear systems on low dimensional nilpotent and solvable Lie groups}, accepted for publication in JDCS, DOI 10.1007/s10883-014-9258-z.
\bibitem{Grushin} V.V.~Grushin \textit{A certain class of hypoelliptic operators}, Mat. Sb.(N.S.)83 (125) (1970) 456-473. 
\bibitem{Jou09} Ph.~Jouan \textit{Equivalence of Control Systems with Linear Systems on
  Lie Groups and Homogeneous Spaces} ESAIM: Control Optimization and Calculus of Variations, 16 (2010) 956-973.
\bibitem{Jou11} Ph.~Jouan \textit{Controllability of linear system on Lie groups}, Journal of Dynamical and control systems, Vol. 17, No 4 (2011) 591-616.
\bibitem{Jou12} Ph.~Jouan \textit{Invariant measures and controllability of finite systems on compact manifolds}, ESAIM: COCV 18 (2012) 643-655.
\bibitem{Jurdjevic97} V.~Jurdjevic \textit{Geometric control theory}, Cambridge university press, 1997.
\bibitem{KL16} V.~Kivioja, E.~Le Donne \textit{Isometries of nilpotent metric groups}, arxiv.org/abs/1601.08172.
\bibitem{Milnor76} J.~Milnor \textit{Curvatures of left invariant metrics on Lie groups}, Advances in Math. 21 (1976), no. 3, 293-329.
\bibitem{Sachkov09} Yu.~L.~Sachkov \textit{Control Theory on Lie groups}, Journal of Mathematical Sciences, Vol. 156, No. 3, 2009.
\bibitem{Takasu} T.~Takasu \textit{Generalized Riemannian Geometry I}, The journal of the Yokohama Municipal University. Series D (1956).
\end{thebibliography}
\end{document}